\documentclass[12pt]{article}%2.12.2017
\usepackage{amsmath}
\usepackage{amssymb}
\usepackage{dsfont}
\usepackage{cite}
\newsymbol \blackbox 1004
\newcommand{\eh}{\hfill}\newlength{\sperr}

\newenvironment{proof}{{\settowidth{\sperr}{\bf\rm
Proof}%
\par\addvspace{0.3cm}\noindent\parbox[t]{1.3\sperr}
{\bf\rm P\eh r\eh o\eh o\eh f\eh }%
}}{\nopagebreak\mbox{}
$\blackbox$\par\addvspace{0.3cm}}

\def\nn{\nonumber}

\def\vk{\varkappa}

\def\Lam{\Lambda}
\def\s{\sigma}

\def\om{\omega}

\def\vp{\varphi}
\def\vt{\vartheta}

\def\wh{\widehat}

\def\ov{\overline}

\def\BC{{\mathbb C}}
\def\BR{{\mathbb R}}
\def\BN{{\mathbb N}}

\def\cla{{\mathcal A}}
\def\clb{{\mathcal B}}
\def\clc{{\mathcal C}}

\def\cll{{\mathcal L}}

\def\im{{\rm Im\ }}

\newcommand{\E}{\mathrm{e}}
\newcommand{\I}{\mathrm{i}}

\newtheorem{Pa}{Paper}[section]
\newtheorem{Tm}[Pa]{{\bf Theorem}}

\newtheorem{Cy}[Pa]{{\bf Corollary}}
\newtheorem{Rk}[Pa]{{\bf Remark}}

\newtheorem{Dn}[Pa]{{\bf Definition}}

\newtheorem{Pn}[Pa]{{\bf Proposition}}

\newenvironment{dedication}
        {\vspace{1ex}\begin{quotation}\begin{center}\begin{em}}    
        {\par\end{em}\end{center}\end{quotation}}

\title{Scattering for general-type Dirac systems \\ on the semi-axis: \\ reflection coefficients and Weyl functions}

\author{Alexander Sakhnovich}

\date{}
\begin{document}
\maketitle

\begin{dedication}
\end{dedication}

\begin{abstract}    We show that for general-type self-adjoint and skew-self-adjoint Dirac systems on the semi-axis
Weyl functions are unique analytic extensions of the reflection
coefficients. New results on the extension of the Weyl functions to the real axis and
on the existence (in the skew-self-adjoint case) of the Weyl functions follow.
Important procedures
to recover general-type Dirac systems from the Weyl functions  are applied  to the recovery of Dirac systems from the reflection coefficients.
We explicitly recover Dirac systems from the rational reflection coefficients as well.

\end{abstract}


\vspace{0.2em}


\section{Introduction}\label{Intro}
\setcounter{equation}{0}
The scattering problems on the semi-axis are of essential interest in theory and
applications (see, e.g., \cite{AHM, AGKS, BeM, GKS6, HryM, Isk, Ism, KoT, Kr, Kur, Mar, Med, SteK1, SteK2, Wed, Zol}).
At the same time, Weyl--Titchmarsh theory on the semi-axis and finite intervals has been actively studied 
(see the books \cite{Tit, Lev, Marb, SaSaR}, recent papers \cite{BeGe, BeKoT, ALS-JST, Zem} and numerous
references therein).
In particular, scattering  and Weyl--Titchmarsh problems for  Dirac systems
\begin{align} &       \label{1.1}
y^{\prime}(x, z )=\I \big(z j+jV(x)\big)y(x,
z ), \quad
x \geq 0 ,
\end{align} 
where $y^{\prime}:=\frac{d}{d x}y$,
\begin{align} &   \label{1.2}
j = \left[
\begin{array}{cc}
I_{m_1} & 0 \\ 0 & -I_{m_2}
\end{array}
\right], \hspace{1em} V= \left[\begin{array}{cc}
0&v\\  \breve v&0\end{array}\right]  \quad ( m_1+m_2=:m),
 \end{align} 
and $m_1,m_2 \in \BN$, are important. 
Here $I_{m_i}$ is the $m_i \times
m_i$ identity
matrix,  $\BN$ stands for the set of natural numbers, $v(x)$ is an $m_1 \times m_2$ matrix function and
$\breve v(x)$ is an $m_2 \times m_1$ matrix function.

The most important cases are the cases of the self-adjoint Dirac systems
\begin{align} &       \label{1.3}
y^{\prime}(x, z )=\I \big(z j+jV(x)\big)y(x,
z ), \quad \breve v(x)  =v(x)^*
\quad (x \geq 0) ,
\end{align} 
and  of the skew-self-adjoint Dirac systems
\begin{align} &       \label{1.4}
y^{\prime}(x, z )=\I \big(z j+jV(x)\big)y(x,
z ), \quad \breve v(x)  =-v(x)^*
\quad (x \geq 0) .
\end{align} 
The Weyl--Titchmarsh theory of the self-adjoint Dirac systems is well-studied (see the references above). It is
also known that Weyl--Titchmarsh (or simply Weyl) functions of the self-adjoint Dirac systems \eqref{1.3} on the semi-axis are closely
related to the scattering data. See, for instance, simple formulas connecting {\it rational} Weyl functions and reflection coefficients
of systems \eqref{1.3}  (where $m_1=m_2$) in \cite{GKS6}, or some special cases of the  scalar system \eqref{1.3} in \cite{BeM}.
Let us mention also the case of the scalar self-adjoint Schr\"odinger equation (see \cite[(1.3)]{BeM} and references in \cite{BeM}).

For the Weyl--Titchmarsh theory of the skew-self-adjoint Dirac systems \eqref{1.4} we refer to the works \cite{ALS90, ClGe, FKRS-LAA}
and \cite[Ch. 3]{SaSaR} (see also references therein). A study of  the Weyl--Titchmarsh theory for systems \eqref{1.4} started much later than for systems \eqref{1.3}
and, in particular,  less is known about the interconnections with the scattering theory.

Therefore, general results on the interconnections between Weyl functions and reflection coefficients, which we obtain and use in this work,
are of interest for both scattering and Weyl--Titchmarsh theories. The next section presents short preliminaries on the scattering theory
of Dirac systems on the semi-axis, and the main section of the paper is Section \ref{RCWF}, where Subsection \ref{SaD} is dedicated
to the self-adjoint systems and Subsection \ref{SkD} is dedicated
to the skew-self-adjoint systems.

Under condition that the entries of $V$ are summable $($integrable$)$ on $\BR_+=[0, \infty)$, that is,
they belong $L^1(\BR_+)$ or, in other words,
\begin{align} &       \label{V}
V(x)\in L_{m \times m}^1(\BR_+),
\end{align} 
we show that for systems \eqref{1.3} and \eqref{1.4} {\it Weyl functions are unique analytic extensions of the reflection
coefficients} $R_L$. (See Theorems \ref{Tm} and \ref{Tmsk}, respectively.) Some new results on the extension of the Weyl functions
of systems \eqref{1.3} on $\BR_+$ (Corollary \ref{Cysa}) and on existence of the Weyl functions of systems \eqref{1.4} are immediate. 
Moreover, important procedures
to recover Dirac systems from the Weyl functions may be used to recover Dirac systems from the reflection coefficients
(see Remarks \ref{RkInvPrsa} and \ref{Rksk}). We {\it explicitly} recover Dirac systems from rational reflection coefficients as well.
For self-adjoint Dirac systems it is the case where \eqref{V} does not necessarily holds but $V(x)\in L_{m \times m}^2(\BR_+)$.
The explicit construction of the Jost solution $F_L$ is new in the skew-self-adjoint case even if $m_1=m_2=1$. 
In view of the recent works on the Weyl functions for systems on the whole axis (see, e.g., \cite{BEKT, GeZ}), the results could be useful also for the scattering
on the axis.

{\it When other requirements are not specified, we  assume  that \eqref{V} holds.}

 In the paper, $\BR$ denotes the real axis, $\BR_+$ denotes the semi-axis $[0, \infty)$, $\BC$ stands for the complex plane,
$\BC_+$ ($\BC_-$) stands for the open upper (lower) half-plane, and $\ov {\BC_+}$ ($\ov {\BC_-}$) stands for the closed upper (lower)
half-plane.
The class of functions, which are summable on $\BR_+$, is denoted by $L^1(\BR_+)$, and the class of functions, which
are square-integrable on each finite interval of $\BR_+$, is denoted by $L^2_{loc}(\BR_+)$.
The notation $\s(A)$ stands for the spectrum of the matrix $A$ and Im$(A)$ denotes the image of $A$.

%The notation $\s(\a)$ stands for the spectrum of the matrix $\a$ and the notation $\diag\{d_1, d_2, \ldots\}$
%stands for the block diagonal matrix with the blocks $d_1, d_2, \ldots$ on the main diagonal.

\section{Preliminaries: Jost solution \\ and reflection coefficient}\label{Prel}
\setcounter{equation}{0}
 Scattering problems on the axis have been actively studied both in this and in the previous centuries
(see, e.g., various important references in \cite{Abl, AHM2, BDT, BDZ, Bond, BRyS, ReeS, SaA8}).
In particular, when the entries of $V(x)$  belong $L^1(\BR)$, there is a plethora
of results concerning Jost solutions and reflection coefficients of Dirac system (on the axis), which hold for 
the more studied case $m_1=m_2$ as well as for the case where $m_1$ is not necessarily equal $m_2$.
The corresponding results, which we discuss in this section, are conveniently collected, for instance,
in \cite[Ch. 3]{Dem}. We formulate those results in terms of Dirac system \eqref{1.1} on the semi-axis.
(Indeed, since $V(x)\in L_{m \times m}^1(\BR_+)$ in our situation, $V(x)$  may be extended to the whole line
by the equality $V(x)\equiv 0$ for $x<0$, and the results for the semi-axis are immediate from the
results on the axis.)

It is known that there is a unique fundamental solution $F_{L  }$ of the system \eqref{1.1} with $z \in \BR$ (Jost solution), which satisfies the asymptotic
relation
\begin{align} &       \label{2.1}
F_{L  }(x,z)=\E^{\I x z j}\big(I_m+o(1)\big), \quad x \to \infty.
\end{align} 
Moreover, the matrix function
\begin{align} &       \label{2.2}
Y(x,z)=\begin{bmatrix}Y_{1}(x,z) \\ Y_2(x,z)  \end{bmatrix}:=F_{L  }(x,z) \begin{bmatrix}I_{m_1} \\ 0  \end{bmatrix}
\end{align} 
can be extended onto $\BC_+$ so that it remains (for each fixed $z\in \BC_+$) a matrix solution of \eqref{1.1}, is  continuous with respect to  $z\in \ov{\BC_+}$
and is analytic on $\BC_+$. This extension is also denoted by $Y(x,z)$.

The norm $\|\E^{-\I x z}Y(x,z)\|$ is uniformly bounded (for $x \in \BR_+$ and  $z \in \ov{\BC_+}$), and the relation
\begin{align} &       \label{2.3}
\lim_{x\to \infty}\E^{-\I x z}Y(x,z)=\begin{bmatrix}I_{m_1} \\ 0  \end{bmatrix} \quad (z \in \ov{\BC_+})
\end{align} 
holds.
We have also
\begin{align} &       \label{2.4-}
\lim_{z \to \infty}\E^{-\I xz}Y_1(x,z)=I_{m_1} \quad {\mathrm{for}} \quad  z\in \ov{\BC_+}.
\end{align}
The reflection coefficient $R_L(z)$ (more precisely, the left reflection coefficient; see, e.g., \cite{Dem,GKS6}) is defined, in the points of invertibility of $Y_1(0,z)$
on $\BR$,
by the equality
\begin{align} &       \label{2.5}
R_L(z)=Y_2(0,z)Y_1(0,z)^{-1}.
\end{align} 
Finally, the following inequality holds in the symmetric case:
\begin{align} &       \label{2.4}
\det\big(Y_1(0,z)\big)\not=0 \quad {\mathrm{if}} \quad  z\in \BR, \,\, \breve v=v^*.
\end{align} 
%%%%%%%%%%%%%%%%%%%%%%%%%%%%%%%%%%%%%%%%%%%%%%%%%%%%%%%%%%%%%%
\section{Reflection coefficients and Weyl functions}\label{RCWF}
\setcounter{equation}{0}
\subsection{Self-adjoint Dirac system}\label{SaD}
\paragraph{1.}
Consider self-adjoint system \eqref{1.3} and suppose first that the entries of $V$ are  locally summable.
The notation $u(x,z)$ stands in this subsection for the fundamental solution of \eqref{1.3}
normalized by the condition
\begin{align} &   \label{3.1}
u(0,z)=I_m.
\end{align}
\begin{Dn} \label{defWeyl} A Weyl--Titchmarsh $($or simply Weyl$)$ function of Dirac system  \eqref{1.3} on $[0, \, \infty)$,
where  the potential $V$ is locally summable, is an  $m_2\times~m_1$ matrix function $\vp$
which satisfies the inequality
\begin{align} &      \label{3.2}
\int_0^{\infty}
\begin{bmatrix}
I_{m_1} & \vp (z)^*
\end{bmatrix}
u(x,z)^*u(x,z)
\begin{bmatrix}
I_{m_1} \\ \vp (z)
\end{bmatrix}dx< \infty , \quad z\in \BC_+.
\end{align} 
\end{Dn}
\begin{Rk}\label{RkWT} In \cite[Subsection 2.2.1]{SaSaR} $($see also \cite{FKRS}$)$, we show that the Weyl function $\vp (z)$ described in Definition \ref{defWeyl} always exists on $\BC_+$, and that it is unique,
holomorphic and contractive $($i.e., $\vp (z)^*\vp (z)\leq I_{m_1})$.
\end{Rk}
\begin{Rk}\label{m1m2} 
In the case $m_1=m_2=:p$, it is convenient to consider fundamental solution $u(x,z)$ normalized by the condition $u(0,z)=K$,
where $K=\begin{bmatrix} I_p & - \I  I_p \\ -I_p & - \I  I_p \end{bmatrix}$.  Then, the Weyl function $\vp=\vp_K$ determined by
\eqref{3.2} is a Herglotz function on $\BC_+$ instead of being contractive. Moreover, the weight function in the Herglotz representation of $\vp_K$ is the spectral
function of the Dirac system \eqref{1.3} $($see, e.g., \cite{ALS-Dir, SaSaR}$)$. Clearly, $\vp$ from Definition \ref{defWeyl} and $\vp_K$
are connected by a simple linear fractional transformation.
\end{Rk}
When the Weyl functions $\vp_K$,
in the sense of Remark \ref{m1m2},  are rational, it is shown in \cite{GKS6} that $R_L(z)=-(I_p+\I \vp_K(z))(I_p-\I \vp_K(z))^{-1}$.
Here, we do not require that  the Weyl function of Dirac system is rational  (and that $m_1=m_2$).
Thus, we prove a much more general statement. The following theorem holds for the Weyl functions $\vp$
described in Definition \ref{defWeyl}.
\begin{Tm}\label{Tm}
Assume that Dirac system \eqref{1.3} satisfies \eqref{V}. Then, its reflection coefficient $R_L(z)$ exists on $\BR$ and
may be continuously extended onto $\BC_+$ so that $R_L(z)$ is analytic on $\BC_+$. Moreover, this extension is unique and the
equality $R_L(z)=\vp(z)$, where $\vp$ is the Weyl function, holds. That is, the unique analytic extension of the reflection coefficient 
coincides on $\BC_+$ with the Weyl function of the Dirac system.
\end{Tm}
\begin{proof}. Recall that, according to Section \ref{Prel}, the matrix function $Y(x,z)$ (given by \eqref{2.2} for the case $z\in \BR$)  is a matrix solution of \eqref{1.3}, 
which is  continuous with respect to   $z\in \ov{\BC_+}$ and is analytic on $\BC_+$. From Section \ref{Prel} we know also that the reflection coefficient exists and
is well-defined by the equality $R_L(z)=Y_2(0,z)Y_1(0,z)^{-1}$ for $z\in \BR$ (see \eqref{2.5} and \eqref{2.4}). 
Formula \eqref{2.4} yields
\begin{align} &      \label{3.3'}
\det\big(Y_1(0,z)\big)\not=0 \quad {\mathrm{for}} \quad z\in \BC_+ ,
\end{align} 
excluding, possibly, some isolated points.
Therefore, the same formula 
$$R_L(z)=Y_2(0,z)Y_1(0,z)^{-1}$$
 continuously  extends  $R_L(z)$ onto $\BC_+$. Moreover, this extension is meromorphic in $\BC_+$
 (and we denote it by $R_L(z)$ similar to the reflection coefficient on $\BR$). Note that (in view of the well-known Luzin--Privalov theorem
\cite{Priv})  the meromorphic extensions of $Y$ and $R_L$  are unique. Finally, if $R_L(z)=\vp(z)$ in $\BC_+$, then $R_L(z)$ is not only meromorphic
but analytic. Thus, it remains to prove the equality
\begin{align} &      \label{!}
Y_2(0,z)Y_1(0,z)^{-1}=\vp(z) \quad (z\in \BC_+).
\end{align} 

From \eqref{2.3} we have 
\begin{align} &      \label{3.3}
Y(x,z)\in L^2_{m \times m_1}(\BR_+)
\end{align} 
for each $z \in \BC_+$. 
On the other hand, formula (2.43) from \cite {SaSaR} implies that the subspace $\cll(z)$ of the span  of columns of $u(x,z)$ ($z\in \BC_+$), 
which belongs to $L^2_m(\BR_+)$,
has dimension less or equal to $m_1$. Thus, taking into account \eqref{3.2} and Remark \ref{RkWT} we see that
\begin{align} &      \label{3.4}
\cll(z) = \im\left(u(x,z)
\begin{bmatrix}
I_{m_1} \\ \vp (z)
\end{bmatrix}\right), \quad \dim \cll(z)=m_1
\end{align} 
for each $z \in \BC_+$, where $\im$ stands for image.

From \eqref{3.3} and \eqref{3.4}, using analyticity of $Y(x,z)$, $u(x,z)$ and $\vp(z)$, we derive the equality
\begin{align} &      \label{3.5}
Y(x,z)=u(x,z)
\begin{bmatrix}
I_{m_1} \\ \vp (z)
\end{bmatrix}c(z)
\end{align} 
for some analytic on $\BC_+$ square matrix function $c(z)$. Relations \eqref{3.1}, \eqref{3.3'} and \eqref{3.5}
imply that $\det c(z)\not=0$ excluding, possibly, some isolated points and that \eqref{!} is valid.
\end{proof}
The next corollary is immediate.
\begin{Cy} \label{Cysa} Let Dirac system \eqref{1.3} satisfy \eqref{V}. Then, the Weyl function of this system
may be continuously extended from $\BC_+$  to the boundary $\BR$.
\end{Cy}
\begin{Rk}\label{RkInvPrsa}  When the entries of $V$ belong to $L^1(\BR_+)\cap L^2_{loc}(\BR_+)$
system \eqref{1.3} is $($in view of  Theorem \ref{Tm}$)$ uniquely recovered from the reflection coefficient $R_L$ using the procedures of \cite[Theorem 4.4]{ALS-JST}
or of \cite[Corollary 3.4]{ALS-HLV}. $($Both procedures are closely related but  the last steps in those procedures differ.$)$
\end{Rk}

\paragraph{2.} In this paragraph, we consider the case of strictly proper rational reflection coefficients $R_L(z)$. Strictly proper rational $m_2 \times m_1$
matrix functions always admit representation (so called {\it realization}):
\begin{equation}
\label{3.6} R_L(z)=\clc(z I_n- \cla)^{-1}\clb,
\end{equation}
where $\cla$ is an $n \times n$ matrix, $\clc$ is an $m_2 \times n$ matrix and $\clb$ is an $n\times m_1$ matrix.
Further in the text we assume that the realization  \eqref{3.6} is a {\it minimal realization}, that is, the value of $n$
in \eqref{3.6} is minimal (among the corresponding values in different realizations of $R_L$).

The inverse problem to recover $V$ from $R_L$ is solved in this case explicitly using the corresponding
results on Weyl functions and fundamental solutions from \cite{FKRS2013,  ALSstab}. The result and approach are close
to \cite[Theorem 9.4]{GKS6} on the $m_1=m_2$ case. However, we do not assume $m_1=m_2$ and the theorem below contains
also some additional  (to those in \cite[Theorem 9.4]{GKS6}) statements. The following theorem is of interest in our considerations
because {\it the entries of $V$  are not necessarily summable on $\BR_+$ under conditions of that theorem.}

\begin{Tm}\label{TmExplSA} Let $R_L(z)$ be an $m_2\times m_1$ strictly proper rational  and contractive on $\BR_+$ matrix function,
which has no poles  on $\BC_+$. Then $R_L(z)$ is a reflection coefficient of a  Dirac system \eqref{1.3} such that its potential $V$ is recovered in the following way.

First choose a minimal
realization \eqref{3.6} of $R_L$. Then, there is
 a strictly positive solution $X>0$ of the Riccati
equation
\begin{equation}\label{3.7}
X\clc^*\clc X-\I(\cla X -X \cla^*)+\clb \clb^*=0.
\end{equation}
We put
\begin{equation}\label{3.8}
A=\cla+\I \clb \clb^* X^{-1}, \quad S(0)=X, \quad \vt_1=\clb, \quad \vt_2=-\I X \clc^*.
\end{equation}
Now, the matrix function $v$, which determines the potential $V=\begin{bmatrix}0 & v \\ v^* & 0\end{bmatrix} $,
is given by the formula
\begin{align} \label{3.9}&
v(x)=-2\I \vt _{1}^{\, *}\E^{\I x A^{*}} S
(x)^{-1}\E^{\I x
A }\vt_2,
\end{align}
where
\begin{align}   \label{3.10}&
S(x)=S(0)+ \int_{0}^{x} \Lambda(t) \Lambda (t)^{*}dt,
\quad
\Lambda (x)= \begin{bmatrix}  \E^{- \I x A } \vt_{1} 
& \E^{\I x
A } \vt_2 \end{bmatrix}.
\end{align}

Moreover, $v$ given above satisfies relations
\begin{align} \label{3.11}&
v(x)\in L^2_{m_1\times m_2}(\BR_+), \quad v(x) \to 0 \quad {\mathrm{for}} \quad x \to \infty,
\end{align}
and the Weyl function $\vp(z)$ of the corresponding Dirac system coincides with $R_L(z)$ on $\BC_+$.
\end{Tm}
\begin{proof}. Consider realization \eqref{3.6}.
According to \cite[Theorem 3.4]{FKRS2013}, $v$ given by \eqref{3.9} is well-defined and
$\vp(z)=\clc(z I_n- \cla)^{-1}\clb$ is the Weyl function of the Dirac system
determined by \eqref{3.9}.

Next, in order to show that  \eqref{3.6} holds for the reflection coefficient
of the constructed system, we note  that $v$ recovered from the Weyl function is unique
(see, e.g., \cite{ALS-JST} and references therein) but the choice of $A$, $S(0)$, $\vk_1$
and $\vk_2$ in \eqref{3.9} and \eqref{3.10} for this $v$ is not unique.
Taking into account \cite[Remark 3.5]{ALSstab}, we may choose $X$ in \eqref{3.7}
so that the additional condition
\begin{align} \label{3.12}&
\s(A)\subset \ov{\BC_-}, \quad {\mathrm{where}} \quad \s \quad {\mathrm{means \,\, spectrum,}}
\end{align}
is fulfilled. The {\it controllability} of the pair $\{A, \vt_1\}$ follows from \eqref{3.8} and from the minimality
of the realization \eqref{3.6}. That is, we have:
\begin{align} \label{3.13}&
{\mathrm{span}}\bigcup_{k=0}^{n-1}\im (A^k \vt_1)=\BC^n.
\end{align}
According to \cite{ALS94}, a fundamental solution of the system \eqref{1.3}, where $v$ is determined
by \eqref{3.9} is given (for $z\not\in \s(A)$) by the equality
\begin{align} \label{3.14}&
\wh u(x, z)=w_{A}(x, z)\E^{\I x z j}, \quad w_{A}(x, z):=I_m-\I j \Lam(x)^*S(x)^{-1}(A-z I_n)^{-1}\Lam(x).
\end{align}
Here $w_A$ has the form of the Lev Sakhnovich transfer matrix function \cite{SaL1, SaSaR},
and $\Lam$ is given in \eqref{3.10}. (We note that $\wh u$ is not normalized by \eqref{3.1}.)

Using \eqref{3.14}, we will construct the Jost solution $F_{L  }$. First, similar to \cite{GKS6}, we rewrite  $w_A$ as
\begin{align}\label{3.15}
& w_{A}(x, z )
\\ \nn &=
\begin{bmatrix}
I_{m_1}-\I \vt_1^*Q(x)^{-1}(A- zI_n)^{-1}\vt_1 &
-\I \vt_1^*\E^{2\I x A^*}R(x)^{-1}(A-z I_n)^{-1}\vt_2 
\\
\I\vt_2^*\E^{-2\I x A^*}Q(x)^{-1} (A-zI_n )^{-1}\vt_1 &
I_{m_2}+\I \vt_2^*R(x)^{-1}(A-z I_n )^{-1}\vt_2
\end{bmatrix},
\end{align}
where
\begin{equation}\label{3.16}
R(x):=\E^{-\I x A } S(x)\E^{\I x A^{*} },  \quad 
Q(x):=\E^{\I x A } S(x)\E^{-\I x A^{*} }.
\end{equation}
From \cite[Lemma 3.7]{ALSstab}, we see that under conditions \eqref{3.12} and \eqref{3.13} the equalities
\begin{equation}\label{3.17}
\lim_{x \to \infty}Q(x)^{-1}=0, \quad \lim_{x \to \infty}Q(x)^{-1}\E^{2\I xA}\vt_2=0
\end{equation}
hold. In a similar to the proof of \cite[Lemma 3.7]{ALSstab} way we will show that 
there are limits:
\begin{equation}\label{3.18}
\lim_{x \to \infty}R(x)^{-1}=\vk_R, \quad \lim_{x \to \infty}R(x)^{-1}\E^{-2\I xA}\vt_1=0.
\end{equation}
Indeed, it easily follows from \eqref{3.7} and \eqref{3.8} that 
\begin{equation}\label{3.18'}
AS(0)-S(0)A^*=\I \Lam(0)j\Lam(0)^*.
\end{equation}
Hence, the identity
\begin{equation}\label{3.19}
AS(x)-S(x)A^*=\I \Lam(x)j\Lam(x)^*
\end{equation}
follows from \eqref{3.10} (see \cite[(3.6)]{FKRS2013}). In view of \eqref{3.10}, \eqref{3.16} and \eqref{3.19},
we have
\begin{equation}\label{3.20}
R^{\prime}(x)=\E^{-\I x A } \big(\Lam(x)j\Lam(x)^*+\Lam(x)\Lam(x)^*\big)\E^{\I x A^{*} }\geq 0.
\end{equation}
Since $R(x)>0$ and $R^{\prime}(x)\geq 0$, the first limit in \eqref{3.18} exists. 

Relations \eqref{3.10}, \eqref{3.16} and \eqref{3.19}
imply also that
\begin{equation}\label{3.21}
A R(x)-R(x) A^*= 
= 
\I (\E^{-2\I x A } \vt_{1} \vt_{1}^{\, *}\E^{2\I x A^{*}}-
\vt_{2} \vt_{2}^{\, *}).
\end{equation}
Multiplying both parts of \eqref{3.21} by $R^{-1}$ from the left and from the right
and taking into account the first limit in \eqref{3.18}, we obtain
\begin{equation}\label{3.22}
 \I(A^*\vk_R -\vk_R A )+\vk_R \vt_{2} \vt_{2}^{\, *}\vk_R
= \lim_{x \to \infty}
\big(R(x)^{-1}\E^{-2\I x A } \vt_{1} \vt_{1}^{\, *}\E^{2\I x A^{*}}R(x)^{-1}\big).
\end{equation}
On the other hand, formula \eqref{3.20} and the second equality in \eqref{3.10}  
yield
\begin{equation}\label{3.23}
R(x)^{-1}\E^{-2\I x A } \vt_{1} \vt_{1}^{\, *}\E^{2\I x A^{*}}R(x)^{-1}=-\big(R(x)^{-1}\big)^{\prime},
\end{equation}
and so $R(x)^{-1}\E^{-2\I x A } \vt_{1} \vt_{1}^{\, *}\E^{2\I x A^{*}}R(x)^{-1}$ is summable on $\BR_+$. More precisely, using \eqref{3.23}
and the first equality in \eqref{3.18} we obtain
\begin{equation}\label{3.24}
\int_0^{\infty}R(x)^{-1}\E^{-2\I x A } \vt_{1} \vt_{1}^{\, *}\E^{2\I x A^{*}}R(x)^{-1}dx=S(0)^{-1}-\vk_R.
\end{equation}
Since, according to \eqref{3.22} and \eqref{3.24},
the left-hand side of \eqref{3.23} is summable on $\BR_+$ and has a limit when $x$ tends to infinity,
this limit should be zero. That is, the second equality in \eqref{3.18} is also valid.
Now, equalities \eqref{3.15}, \eqref{3.17} and \eqref{3.18} imply that
\begin{align}\label{3.25}&
w_{A}(x, z )
=
\begin{bmatrix}
I_{m_1} &
0
\\
0 & \chi(z)
\end{bmatrix} +o(1) \quad {\mathrm{for}} \quad x\to \infty, 
\\ & \label{3.25'} \chi(z):=I_{m_2}+\I \vt_2^*\vk_R(A-z I_n )^{-1}\vt_2.
\end{align}
In view of the definition \eqref{2.1} (of the Jost solution $F_{L  }$) and of the relations \eqref{3.14} and \eqref{3.25},
we derive
\begin{equation}\label{3.26}
F_{L  }(x,z)=w_A(x,z)\E^{\I x z j}\begin{bmatrix}
I_{m_1} &
0
\\
0 & \chi(z)^{-1}
\end{bmatrix} \quad (z=\ov{z}, \quad z\not\in \s(A)).
\end{equation}
Here the invertibility of $\chi(z)$ is immediate from \eqref{3.25'} and from the equality $\I(A^*\vk_R -\vk_R A )+\vk_R \vt_{2} \vt_{2}^{\, *}\vk_R=0$ (see, e.g.
\cite[Appendix]{GKS6}).
The equality $\I(A^*\vk_R -\vk_R A )+\vk_R \vt_{2} \vt_{2}^{\, *}\vk_R=0$ follows in turn from \eqref{3.18} and \eqref{3.22}.

Taking into account definition \eqref{2.2}  and relations \eqref{3.26} and \eqref{3.15}, we see that $Y(0,z)= w_{A}(0, z )\begin{bmatrix} I_{m_1} \\ 0\end{bmatrix}$,
and so
\begin{align}\nn
 Y_2(0,z)Y_1(0,z)^{-1}=&\I\vt_2^*S(0)^{-1} (A-zI_n )^{-1}\vt_1
\\ & \times  \big(I_{m_1}-\I \vt_1^*S(0)^{-1}(A- zI_n)^{-1}\vt_1\big)^{-1}.
\label{3.27}
\end{align}
Using again \cite[Appendix]{GKS6} and taking into account   \eqref{3.8}, we write
\begin{align}& \nn
   \big(I_{m_1}-\I \vt_1^*S(0)^{-1}(A- zI_n)^{-1}\vt_1\big)^{-1}=I_{m_1}+\I \vt_1^*S(0)^{-1}(A^{\times}- zI_n)^{-1}\vt_1, 
   \\ &
   A^{\times}=A -\I \vt_1\vt_1^*S(0)^{-1}=\cla .
\label{3.28}
\end{align}
Substituting \eqref{3.28} into \eqref{3.27}, after easy transformations we have
\begin{align}&
 Y_2(0,z)Y_1(0,z)^{-1}=-\I\vt_2^*S(0)^{-1} (zI_n-\cla )^{-1}\vt_1=\clc (zI_n-\cla )^{-1}\clb .
\label{3.29}
\end{align}
Thus, definition \eqref{2.5} and formula \eqref{3.29} show that  \eqref{3.6} holds for the reflection coefficient $R_L$.

Finally, rewriting $v(x)$ in the form \cite[(3.17)]{ALSstab} and using 
\cite[(3.9) and (3.13)]{ALSstab} we see that \eqref{3.11} is valid.
\end{proof}

\subsection{Skew-self-adjoint Dirac system}\label{SkD}
\paragraph{1.} In this subsection, $u$ stands for the fundamental solution of the skew-self-adjoint system \eqref{1.4} normalized, similar to the self-adjoint
case, by \eqref{3.1}.
Weyl--Titchmarsh theory for the skew-self-adjoint Dirac system \eqref{1.4} and corresponding references are given
in \cite[Ch. 3]{SaSaR}. Let us recall several basic definitions and facts from \cite{SaSaR}. 
Weyl functions are not necessarily considered on $\BC_+$ but on some open half-plane $\BC_M=\{z: \,\, z\in \BC, \,\, \Im(z)>M \geq 0\}$.
An $m_2 \times m_1$ matrix function $\vp$, which  satisfies (on $\BC_M$ for some  $M>0$) the inequality \eqref{3.2},
is called a Weyl function of the  Dirac system  \eqref{1.4}. Under condition
 \begin{align}\label{3.30}
 \|v(x)\| \leq M  \quad \mathrm{for} \,\, x\in [0, \, \infty),
 \end{align}
the Weyl function $\vp$ exists and is unique on $\BC_M$. It is  holomorphic 
 and contractive as well. Moreover, under condition \eqref{3.30} the Weyl function $\vp$  satisfies the inequality
\begin{align}&      \label{3.31}
\sup_{x \leq l, \, z\in \BC_M}\left\| \E^{-\I x z }u(x,z)\begin{bmatrix}
I_{m_1} \\ \vp(z)
\end{bmatrix} \right\|<\infty
\end{align}
on  any finite interval $[0, \, l]$. 

In the case of the {\it locally bounded} potentials $V$, the property \eqref{3.31} holding (together with analyticity in some fixed $\BC_M$) on any finite interval $[0, \, l]$
is the definition of the so called GW-function (generalized Weyl function) $\vp(z)$. We note that there is no more than one GW-function.

\paragraph{2.} Similar to the  self-adjoint Dirac system, in the case of the skew-self-adjoint Dirac system we are interested again  in the behavior
of $\det Y_1(0,z)$.
\begin{Pn} Let Dirac system  \eqref{1.4} be given and let \eqref{V} hold. 
Then, $\det Y_1(0,z)$ is analytic on $\BC_+$ and continuous on $\ov{\BC_+}$. Moreover,  $\det Y_1(0,z)$ satisfies $($for some $M, \, \wh M \, >0)$
the inequalities
\begin{align}&      \label{3.32}
\det Y_1(0,z)\not=0 \quad (\Im (z) \geq M);
\\ &      \label{3.33}
 \det Y_1(0,z)\not=0 \quad   {\mathrm{on \,\, the \,\, rays}} \quad \{z: \,\, z\in \BR, \,\, |z|>\wh M\};
 \\ &      \label{3.34}
\det Y_1(0,z)\not=0 \quad   {\mathrm{almost \,\, everywhere \,\, on}} \quad \BR.
\end{align}
\end{Pn}
\begin{proof}.  Since (according to ``Preliminaries") $Y(0,z)$  is analytic on $\BC_+$ and continuous on $\ov{\BC_+}$, $\det Y_1(0,z)$ has the same properties.
The inequalities \eqref{3.32} and \eqref{3.33}  are immediate from \eqref{2.4-}. Moreover, in view of  \eqref{3.32} and Luzin--Privalov theorem we see that  \eqref{3.34} holds.
\end{proof}
Thus, the reflection coefficient $R_L$ is well-defined via \eqref{2.5} almost everywhere on $\BR$. Now, the next statement easily follows from \eqref{2.3}.
\begin{Tm} \label{Tmsk} Let Dirac system  \eqref{1.4} be given and let \eqref{V} hold. Then, the reflection coefficient $R_L$ is defined $($via \eqref{2.5}$)$
almost everywhere on $\BR$ and admits a unique, continuous on $ \{z: \,\, z\in \BR, \,\, |z|>\wh M\}\cup \BC_+$ and meromorphic on $\BC_+$ extension
also denoted by $R_L$ and given by \eqref{2.5} $(z \in \BC_+)$. This extension coincides for $z\in \BC_M$ with the Weyl function $\vp(z)$ 
and has the property \eqref{3.31},
that is, coincides with the GW-function as well.
\end{Tm}
\begin{proof}. The properties of $R_L(z)$ on $\BR$ were already shown above. Recalling again that $Y(0,z)$ is analytic on $\BC_+$ and continuous on $\ov{\BC_+}$,
and using \eqref{2.5}, \eqref{3.33} and Luzin--Privalov theorem, we see that $R_L$ admits a unique, continuous on 
$$ \{z: \,\, z\in \BR, \,\, |z|>\wh M\}\cup \BC_+$$ 
and meromorphic 
on $\BC_+$ extension. According to \eqref{3.1}, we have the equality 
\begin{align}&      \label{3.35}
u(x,z)\begin{bmatrix}
I_{m_1} \\ Y_2(0,z)Y_1(0,z)^{-1}
\end{bmatrix} =Y(x,z)Y_1(0,z)^{-1}.
\end{align}
In view of \eqref{3.32}, both parts of \eqref{3.35} are well-defined on $\BC_M$.
Taking into account \eqref{2.3}, we derive that the right-hand side of \eqref{3.35}
belongs to $L^2_{m \times m_1}(\BR_+)$ for  all $z\in \BC_M$. Hence, substituting $\vp(z)=Y_2(0,z)Y_1(0,z)^{-1}$
into the left-hand side of \eqref{3.35} we see that \eqref{3.2} holds. In other words, $\vp(z)=Y_2(0,z)Y_1(0,z)^{-1}=R_L(z)$
is the Weyl function of Dirac system on $\BC_M$. Now, the uniform boundedness of  $\|\E^{-\I x z}Y(x,z)\|$, asymptotics \eqref{2.4-}
and relation \eqref{3.35} yield \eqref{3.31}. (In fact, we obtain even stronger inequality.)
\end{proof}
\begin{Rk}\label{Rksk} According to Theorem \ref{Tmsk}, $R_L$ coincides with the Weyl function, and so the procedures in Theorems 3.21, 3.24 and 3.30
from \cite{SaSaR} may be used to recover the potential $V$ of the skew-self-adjoint Dirac system from $R_L$. Note that the condition of boundedness
or local boundedness of $\| V\|$ $($for solving inverse problems$)$ in \cite[Ch. 3]{SaSaR} may be essentially weakened similar to the way it is done for the self-adjoint case
in \cite{ALS-JST}.
\end{Rk}
\paragraph{3.} In the case of rational $R_L$, we again recover the potential $V$ explicitly.
We explicitly construct  also the Jost solution $F_{L  }$, which is new (for skew-self-adjoint Dirac system)
even if $m_1=m_2=1$. Moreover, the explicit expression for $F_{L  }$ provides another way to show
that reflection coefficients and Weyl functions coincide.
\begin{Tm}\label{TmIpes} Let  $R_L(z)$ 
be a strictly proper rational $m_2\times m_1$ matrix function. 
Then $R_L(z)$ is the reflection coefficient  $($and the Weyl function in $\BC_M$ for some $M>0)$ of  the Dirac system \eqref{1.4}
with a bounded potential $V$.
The potential $V$ is uniquely recovered using the following procedure.

Assuming that \eqref{3.6} is a minimal realization of $\vp(z)$ and choosing a positive solution $X>0$ $($which always exists$)$
of  the Riccati equation
\begin{align} \label{Ric}&
X\clc^*\clc X+\I(\cla  X-X\cla^*)-\clb \clb^*=0,
\end{align}
we put
\begin{align}
\label{3.36}&
A=\cla +\I \clb \clb^*X^{-1},  \quad S(0)=X, \quad \vt_1=\clb, \quad \vt_2=\I X \clc^*.
\end{align}
The potential $V$ corresponding to  $R_L$ has the form
\begin{align} \label{3.37}&
V(x)=\begin{bmatrix}0 & v(x) \\ -v(x)^* & 0 \end{bmatrix}, \quad v(x)=-2 \I  \vt _{1}^{\, *}\E^{\I x A^{*}} S
(x)^{-1}\E^{\I x
A }\vt_2,
\end{align}
where
\begin{align} 
 \label{3.38}& 
S(x)=S(0)+ \int_{0}^{x} \Lambda(t) j \Lambda (t)^{*}dt,
\quad
\Lambda (x)= \begin{bmatrix}  \E^{- \I x A } \vt_{1} 
& \E^{\I x
A } \vt_2 \end{bmatrix}.
\end{align}
\end{Tm}
\begin{proof}. The fact that the strictly proper rational matrix function with the minimal realization \eqref{3.6}
is the  Weyl function of the Dirac system \eqref{1.4} where $V$
has the form \eqref{3.37}, \eqref{3.38}, is proved in \cite[Theorem 3.9]{FKKS} (see also \cite[Theorem 2.9]{FKRS-LAA}). 
According to \cite[Corollary 3.6]{FKRS-LAA}, we have $\lim_{x\to \infty}V(x)=0$ (and so $V(x)$ is bounded in the norm).
After taking into account \cite[Corollary 3.25]{SaSaR}, the uniqueness of the solution of the inverse problem
(i.e., the problem to recover $V$) in the class of bounded
on $\BR_+$ potentials follows from \cite[Theorem 3.21]{SaSaR}.

Now, we construct the Jost solution of  the system \eqref{1.4}, \eqref{3.37} and show that $R_L$ of the form \eqref{3.6} is the reflection coefficient of  \eqref{1.4},
\eqref{3.37}.
The matrix function 
\begin{align} \label{3.39}&
\wh u(x, z)=w_{A}(x, z)\E^{\I x z j} \quad (w_{A}(x, z):=I_m-\I \Lam(x)^*S(x)^{-1}(A-z I_n)^{-1}\Lam(x))
\end{align}
is a fundamental solution of \eqref{1.4}, \eqref{3.37} (see \cite[(2.9)]{FKKS}). Similar to \eqref{3.15},
we rewrite the blocks of $w_A$:
\begin{align}\label{3.40}
& w_{A}(x, z )
\begin{bmatrix}
0 \\ I_{m_2}
\end{bmatrix}
=
\begin{bmatrix}
-\I \vt_1^*\E^{2\I x A^*}R(x)^{-1}(A-z I_n)^{-1}\vt_2 
\\
I_{m_2}-\I \vt_2^*R(x)^{-1}(A-z I_n )^{-1}\vt_2
\end{bmatrix}, 
\end{align}
where $R(x):=\E^{-\I x A } S(x)\E^{\I x A^{*} }$.  In view of Proposition 3.5 and equality (3.12) from \cite{FKRS-LAA},
we have
\begin{equation}\label{3.41}
\lim_{x \to \infty}R(x)^{-1}=0, \quad \lim_{x \to \infty}R(x)^{-1}\E^{-2\I xA}\vt_1=0.
\end{equation}
Formulas \eqref{3.40} and \eqref{3.41} yield
\begin{align}\label{3.42}
& \lim_{x\to \infty}\left(w_{A}(x, z )
\begin{bmatrix}
0 \\ I_{m_2}
\end{bmatrix}\right)
=
\begin{bmatrix}
0 \\ I_{m_2}
\end{bmatrix}.
\end{align}
Moreover, $w_A(x,z)$ is unitary for $z\in \BR$ \cite{SaL1} (see also \cite[Corollary 1.13]{SaSaR} or \cite[(3.28)]{FKKS}).
Therefore, formula \eqref{3.42} implies that
\begin{align}\label{3.43}
& \lim_{x\to \infty}\left(\begin{bmatrix}
0 & I_{m_2}
\end{bmatrix}w_{A}(x, z )
\begin{bmatrix}
I_{m_1} \\ 0
\end{bmatrix}\right)
=0 \quad (z\in \BR).
\end{align}
It remains to study the upper left block $\big(w_A\big)_{11}$ of $w_A$. This block admits a representation, which is similar to \eqref{3.15}:
\begin{align}\label{3.44}
&\big(w_A\big)_{11}(x, z )
=
I_{m_1}-\I \vt_1^*Q(x)^{-1}(A- zI_n)^{-1}\vt_1 , \quad Q(x):=\E^{\I x A } S(x)\E^{-\I x A^{*} }.
\end{align}
The matrix identity for $S(x)$ in the skew-selfadjoint case slightly differs from \eqref{3.19} and has the form
\begin{equation}\label{3.45}
AS(x)-S(x)A^*=\I \Lam(x)\Lam(x)^*,
\end{equation}
see \cite[(2.5)]{FKKS}. Hence, taking into account formula \eqref{3.38}, the second equality in \eqref{3.44} and the identity \eqref{3.45},
we derive
\begin{align} & \label{3.46}
Q^{\prime}(x)=\E^{\I x A }\Lam(x)(j-I_m)\Lam(x)^*\E^{-\I x A^{*} }\leq 0, \\
& \label{3.47}
AQ(x)-Q(x)A^*=\I \big(\vt_1\vt_1^*+\E^{2 \I x A }\vt_2\vt_2^*   \E^{-2 \I x A^* }\big).
\end{align}
In view of \eqref{3.46}, $Q(x)$ monotonically decreases: 
\begin{equation}\label{3.48}
\lim_{x\to \infty}Q(x)=Q_{\infty} \geq 0.
\end{equation}
Since the realization \eqref{3.6} is minimal
and $A$, $\vt_1$ are given by \eqref{3.36}, the pair $\{A,\vt_1\}$ is controllable. Thus, we have (see \cite[(2.12)]{FKRS-LAA}):
\begin{equation}\label{3.49}
\s(A)\subset \BC_+.
\end{equation}
Relations \eqref{3.47}--\eqref{3.49} imply that
\begin{equation}\label{3.50}
AQ_{\infty}-Q_{\infty}A^*=\I \vt_1\vt_1^*.
\end{equation}
From \eqref{3.49} and \eqref{3.50} (in a standard way, using residues), we obtain
\begin{equation}\label{3.51}
Q_{\infty}=\frac{1}{2\pi}\int_{-\infty}^{+\infty}(A-t I_n)^{-1} \vt_1\vt_1^*(A^*-t I_n)^{-1}dt.
\end{equation}
Recall that the pair $\{A,\vt_1\}$ is controllable, and so \eqref{3.51} implies $Q_{\infty}>0$.
Therefore, formulas \eqref{3.44} and \eqref{3.48} yield
\begin{align}\label{3.52}
&\lim_{x\to \infty} \big(w_A\big)_{11}(x, z )
=
I_{m_1}-\I \vt_1^*Q_{\infty}^{-1}(A- zI_n)^{-1}\vt_1 :=\om(z).
\end{align}
Taking into account the definition \eqref{2.1} of the Jost solution 
and relations \eqref{3.39}, \eqref{3.42}, \eqref{3.43} and \eqref{3.52},
we write down an explicit expression for the Jost solution:
\begin{align}&
\label{3.53}
F_{L  }(x,z)=w_A(x,z) \E^{\I x z j} \begin{bmatrix} \om(z)^{-1} & 0 \\ 0 & I_{m_2} 
\end{bmatrix}.
\end{align}
It is immediate from \eqref{3.53} that the equality
\begin{align}
\label{3.54}
R_L(z)=\big(w_A\big)_{21}(0,z) \big(w_A\big)_{11}(0,z)^{-1}= &-\I\vt_2^*S(0)^{-1} (A-zI_n )^{-1}\vt_1
\\ \nn & \times  \big(I_{m_1}-\I \vt_1^*S(0)^{-1}(A- zI_n)^{-1}\vt_1\big)^{-1}
\end{align}
holds for the reflection coefficient when $z\in \BR$. Hence, quite similar to the proof of \eqref{3.28} and \eqref{3.29}
one may show that
\begin{align}
\label{3.55}
R_L(z)=\I\vt_2^*S(0)^{-1} (zI_n-\cla )^{-1}\vt_1=\clc (zI_n-\cla )^{-1}\clb.
\end{align}
That is, $R_L$ given by \eqref{3.6} is the reflection coefficient of the system \eqref{1.4} with the potential $V$
given by \eqref{3.37} and \eqref{3.38}.

\end{proof}

\bigskip

\noindent{\bf Acknowledgments.}
 {This research   was supported by the
Austrian Science Fund (FWF) under Grant  No. P29177.}
%%%%%%%%%%%%%%%%%%%%%%%%%%%%%%%%%%%%%%%%%%%%%%
%%%%%%%%%%%%%%%%%%%%%%%%%%%%%%%%%%%%%%%%%%%%%%%

\begin{flushright}

A.L. Sakhnovich,\\
Fakult\"at f\"ur Mathematik, Universit\"at Wien, \\
Oskar-Morgenstern-Platz 1, A-1090 Vienna, Austria\\
e-mail: oleksandr.sakhnovych@univie.ac.at
\end{flushright}

\end{document}